# New Large Graphs with Given Degree and Diameter *

F.Comellas and J.Gómez
Departament de Matemàtica Aplicada i Telemàtica
Universitat Politècnica de Catalunya

## Abstract

In this paper we give graphs with the largest known order for a given degree $\Delta$ and diameter $D$. The graphs are constructed from Moore bipartite graphs by replacement of some vertices by adequate structures. The paper also contains the latest version of the $(\Delta, D)$ table for graphs.

## 1 Introduction

A question of special interest in Graph theory is the construction of graphs with an order as large as possible for a given maximum degree and diameter or $(\Delta, D)$-problem. This problem receives much attention due to its implications in the design of topologies for interconnection networks and other questions such as the data alignment problem and the description of some cryptographic protocols.

The $(\Delta, D)$ problem for undirected graphs has been approached in different ways. It is possible to give bounds to the order of the graphs for a given maximum degree and diameter (see [5]). On the other hand, as the theoretical bounds are difficult to attain, most of the work deals with the construction of graphs, which for this given diameter and maximum degree have a number of vertices as close as possible to the theoretical bounds (see [4] for a review).

Different techniques have been developed depending on the way graphs are generated and their parameters are calculated. Many $(\Delta, D)$

*Research supported in part by the *Comisión Interministerial de Ciencia y Tecnologia*, Spain, under grant TIC90-0712

graphs correspond to Cayley graphs [6, 11, 7] and have been found by computer search. Compounding is another technique that has proved useful and consists of replacing one or more vertices of a given graph with another graph or copies of a graph and rearranging the edges suitably (see [13, 14]). Compound graphs and Cayley graphs constitute most of the large $(\Delta, D)$ graphs described in the literature.

Other large graphs have been found as graph products or special methods. For instance, a graph on an alphabet may be constructed as follows: the vertices are labeled with words on a given alphabet and a rule that relates two different words gives the edges. This representation of the graph is useful for the direct determination of the diameter. As an example, the Lente graph -the best $(5,3)$ graph known- was found as a graph on an alphabet.

An extensive part of the search for large graphs has been performed using computer methods. Usually, the computer is used for generating the graphs and testing for the desired properties. If an exhaustive search is not possible due to the extent of the state space of possible solutions some authors use local search or random algorithms (for example [7]).

In this paper compounding is used to construct new families of large $(\Delta, D)$ graphs that considerably improve the results known for diameter 6. The technique is a generalization of a method used by Quisquater [16], based on the replacement by a complete graph of a single vertex from a bipartite Moore graph. Gómez, Fiol and Serra in [14] modified the technique in order to replace several vertices. This paper extends this technique so that we are able to present a general rule for the replacement of a large number of vertices.

Section 2 is devoted to notation and some known results concerning Moore bipartite graphs. In Section 3 we describe the general technique for the construction of large graphs with diameter 6 and in Section 4 we construct special graphs based also on bipartite Moore graphs. Finally, we present an updated version of the table of large $(\Delta, D)$ graphs.

## 2 Notation and known results

A graph, $G = (V, A)$, consists of a non empty finite set $V$ of elements called *vertices* and a set $A$ of pairs of elements of $V$ called *edges*. The number of vertices $N = |G| = |V|$ is the *order* of the graph. If $(x, y)$ is an edge of $A$, we say that $x$ and $y$ (or $y$ and $x$) are *adjacent* and this is usually written $x \sim y$. It is also said that $x$ and $y$ are the *endvertices* of the edge $(x, y)$. The graph $G$ is *bipartite* if $V = V_1 \cup V_2$ and $V_1 \cap V_2 = \emptyset$ and any edge $(x, y) \in E$ has one endvertex in $V_1$ and the other in $V_2$. The *degree* of a vertex $\delta(x)$ is the number of vertices adjacent to $x$. The *degree* of $G$ is $\Delta = \max_{x \in V} \delta(x)$. A graph is regular of degree $\Delta$ or $\Delta$-*regular* if the degree of all vertices equal $\Delta$. The *distance* between two vertices $x$ and $y$, $d(x, y)$, is the number of edges of a shortest path between $x$ and $y$, and its maximum value over all pair of vertices, $D = \max_{x,y \in V} d(x, y)$, is the *diameter* of the graph. A $(\Delta, D)$ *graph* is a graph with maximum degree $\Delta$ and diameter at most $D$.

The order of a graph with degree $\Delta$ ($\Delta > 2$) of diameter $D$ is easily seen to be bounded by
$$1 + \Delta + \Delta(\Delta - 1) + \ldots + \Delta(\Delta - 1)^{D-1} = \frac{\Delta(\Delta - 1)^D - 2}{\Delta - 2} = N(\Delta, D)$$
This value is called the *Moore bound*, and it is known that, for $D \geq 2$ and $\Delta \geq 3$, this bound is only attained if $D = 2$ and $\Delta = 3, 7$, and (perhaps) 57, (see [5]). In this context, it is of great interest to find graphs which for a given diameter and maximum degree have a number of vertices as close as possible to the Moore bound.

In this paper we propose a way of modifying some known large bipartite graphs in order to obtain new larger graphs for some values of the degree and the diameter.

For bipartite graphs, by counting arguments it is easy to obtain the following upper bound (see [5]) for the maximum order of a $(\Delta, D)$ graph:
$$b(\Delta, D) \leq 2 \frac{(\Delta - 1)^D - 1}{\Delta - 2}, \Delta > 2$$
The bipartite graphs that attain the bound are called *bipartite Moore graphs*. They exist only for $D = 2$ (complete bipartite graphs $K_{\Delta, \Delta}$) or $D = 3, 4, 6$. For these values of $D$, bipartite Moore graphs exist if $\Delta - 1$ is a prime power ([5, 3]). For $D = 4$ they are called *generalized quadrangles*, $Q_q$, and are graphs whose vertices are the points and lines of a non-degenerate quadric surface in projective 4-space with two vertices being adjacent if and only if they correspond to an incident point-line pair in the surface. The graphs $Q_q$ have order $N = 2\frac{q^4 - 1}{q - 1}$ and degree $\Delta = q + 1$. For $D = 6$ the bipartite Moore graphs are called *generalized hexagons*, $H_q$, and they are defined in a similar way to $Q_q$, see [2]. They have order $N = 2\frac{q^6 - 1}{q - 1}$ and degree $\Delta = q + 1$.

Here follow some known results that will be used in the next Sections. If $G$ is a $(\Delta, D)$ bipartite graph and $d(x,y) = D$, $x, y \in V(G)$ then there exist $\Delta$ disjoint paths between $x$ and $y$ of length exactly D. Besides, if $G = (V_1 \cup V_2, E)$ is any bipartite graph of even (odd) diameter $D$, the distance between $x \in V_1$ and any $y \in V_2$ ($y \in V_1$) is at most $D - 1$.

## 3    Large graphs from generalized hexagons

In this section we give new large graphs of diameter 6 obtained from modifications of generalized hexagons $H_q$. We obtain the largest known graphs with diameter 6 and degree $\Delta = q+1$ when $q$ is a prime power. The technique consists of expanding the bipartite Moore graph $H_q$ by replacing several vertices by complete graphs and creating some new adjacencies.

Let us consider the subgraph of $H_q$ shown in Figure 1.

We replace some vertices $x_{i\,j\,k}$ by copies of the complete graph with $h$ elements $K_h$ ($h < \Delta$) that will be denoted $K_{i\,j\,k}$. These replacements must verify the following conditions:

a. Each vertex of a $K_{i\,j\,k}$ graph has at least one edge to one of the vertices that were adjacent to $x_{i\,j\,k}$ (and not shown in Figure 1).

b. Each graph $K_{i\,j\,k}$ has (at least) an edge to $K_{i\,j\,l}$, $l \neq k$.

c. Each graph $K_{i\,j\,k}$ has (at least) an edge to $K_{i\,m\,n}$, $j \neq m$.

d. Each vertex of $K_{i\,j\,k}$ is adjacent, at least, to one vertex of $K_{o\,p\,m}$, $o \neq i$.

e. The different graphs $K_{i\,j\,k}$ are interconnected in such a way that the degree of each vertex is not greater than $\Delta$.

The new graph is denoted $H_q(K_h)$.

**Lemma.** *The graph $H_q(K_h)$ has diameter $D = 6$.*
**Proof.** First, we must observe that after replacement of the vertices $x_{i\,j\,k}$ of $H_q$:

i. A path of maximum length 6, such that it has no vertices $x_{i\,j\,k}$, is unaltered.

ii. The shortest path between any two vertices will increase its length by at most two (the maximum number of vertices $x_{i\,j\,k}$ that might contain).

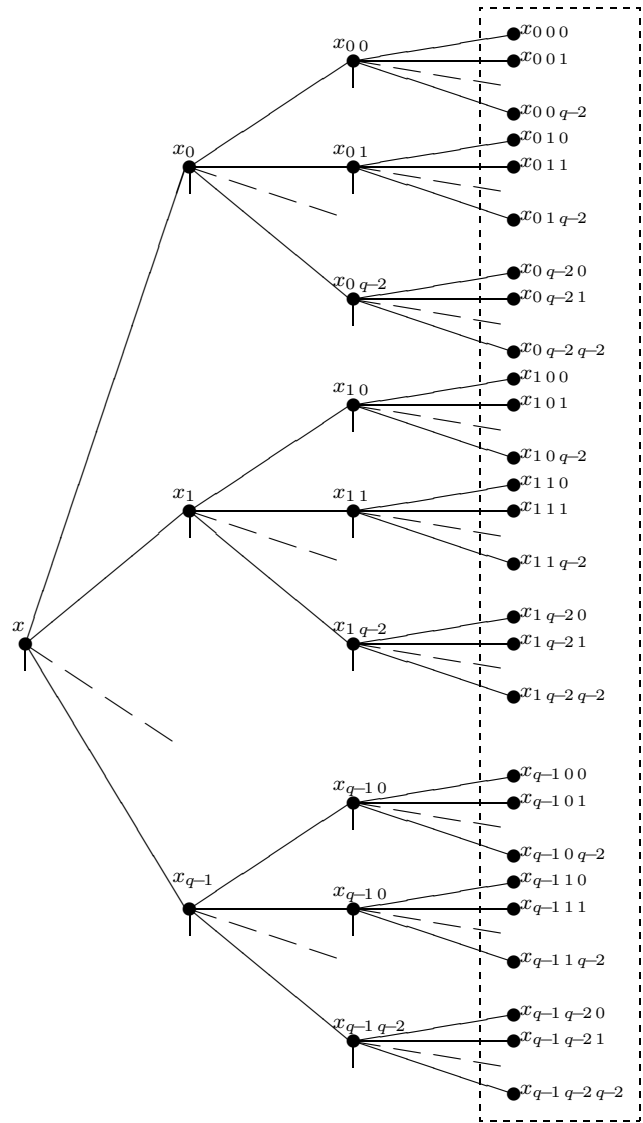

Figure 1: The subgraph of $H_q$ to be modified.

iii. Because of condition (b) the maximum distance between two vertices of $K_{i\,j\,k}$ and $K_{i\,j\,l}$ is 3.

We have the following cases:

1. Let us consider two vertices in $H_q$ at distance 6 and such that one of them, at least, is not replaced. Then, as there exist $\Delta$ disjoint paths of length 6 between these vertices, from condition (a) and Figure 1, we can always find one path among these that is unaffected by the replacements.

2. From conditions (b) and (c), the distance between vertices of $K_{i\,j\,k}$ and $K_{i\,m\,n}$ is at most 5.

3. From conditions (b), (c) and (d), the distance between vertices of $K_{i\,j\,k}$ and $K_{o\,p\,n}$ is at most 6.

4. Let us consider a path of maximum length 5 in $H_q$ between a vertex that will not be modified and any other vertex. If this path contains a vertex $x_{i\,j\,k}$, from (iii) and (2), this path increases its length by at most one unit. $\square$

We have contructed new large $(\Delta, D)$ graphs with diameter 6 and degrees 6,8,9,10,12 and 14 as follows:

$H_5(K_4)$ is obtained by replacing in $H_5$ vertices $x_{i\,j\,k}$, $i = 0$ and $0 \le j, k, \le 3$, by complete graphs $K_4$. In the same way the vertices to be replaced to obtain $H_7(K_6)$ have indices $i = 0$ and $0 \le j, k, \le 5$; $H_8(K_6)$, $0 \le i \le 1$, $0 \le j \le 5$, $0 \le k \le 4$; $H_9(K_6)$, $0 \le i \le 1$, $0 \le j, k, \le 7$; $H_{11}(K_6)$, $0 \le i \le 2$, $0 \le j, k, \le 9$ and $H_{13}(K_7)$, $0 \le i \le 3$, $0 \le j \le 11$, $0 \le k \le 10$.

## 4 Special constructions

In this section we present new large graphs that, as in the previous Section, are obtained from bipartite Moore graph by replacement of some vertices for complete graphs.

First, we modify $Q_4$, the bipartite Moore graph of diameter 4 and degree 5, and without increasing the diameter, we obtain a new large $(5, 4)$ graph.

Let us consider any vertex of $Q_4$ as a starting point for the modification. From this vertex we consider four of its adjacent vertices and for each one of them two adjacent vertices at distance two from the initial vertex. These eight vertices are replaced by complete graphs $K_3$ and new edges are added as shown in Figure 2. As a result, we have created a new graph that we call $Q_4(K_3)$ with 16 vertices more than $Q_4$. Reasoning in a similar way to in Section 3, it is easy to check that the diameter of $Q_4(K_3)$ is 4.

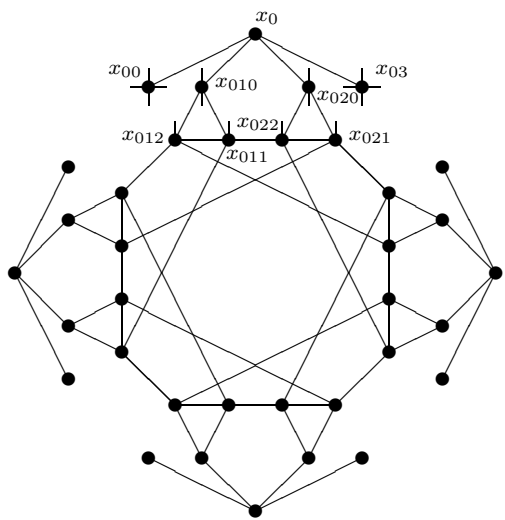

Figure 2: Graph $Q_4(K_3)$ with 186 vertices.

Two new large graphs with diameter 6 and degrees 4 and 5 may be constructed following the technique described in Section 3 but without using condition (d). This is a sufficient condition to ensure that vertices of the different complete graphs are at a distance no greater than 6. Because we do not use this condition, we are able to substitute more vertices, but we must join the complete graphs by a special arrangement of edges such that the diameter of the resulting graph is still not greater than 6. More precisely, Figure 3 shows the modification of $H_3$ that replaces 6 vertices by 6 copies of $K_3$ giving $H_3(K_3)$ with 740 vertices and Figure 4 shows how to obtain from $H_4$ a graph, $H_4(K_4)$, with 2754 vertices.

Finally, a large $(7,6)$ graph was obtained from $H_4(K_4)$ and $H_5$ using a compound technique described in [14].

These four graphs, together with the improvements obtained in Section 3, are displayed in boldtype in Table 1. Table 2 contains the description of the entries of Table 1. An updated version of these tables

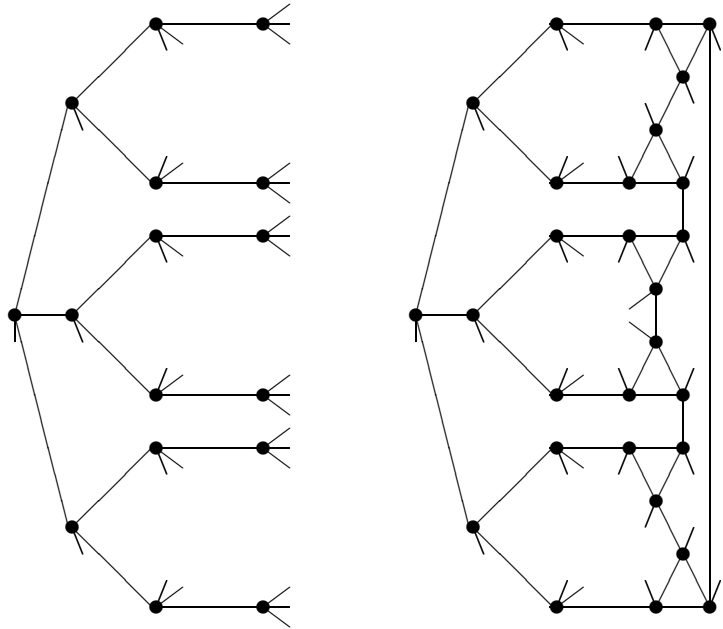

Figure 3: Modification of $H_3$ that gives $H_3(K_3)$

is available from J-C Bermond and C. Delorme (e-mail: cd@lri.fr).

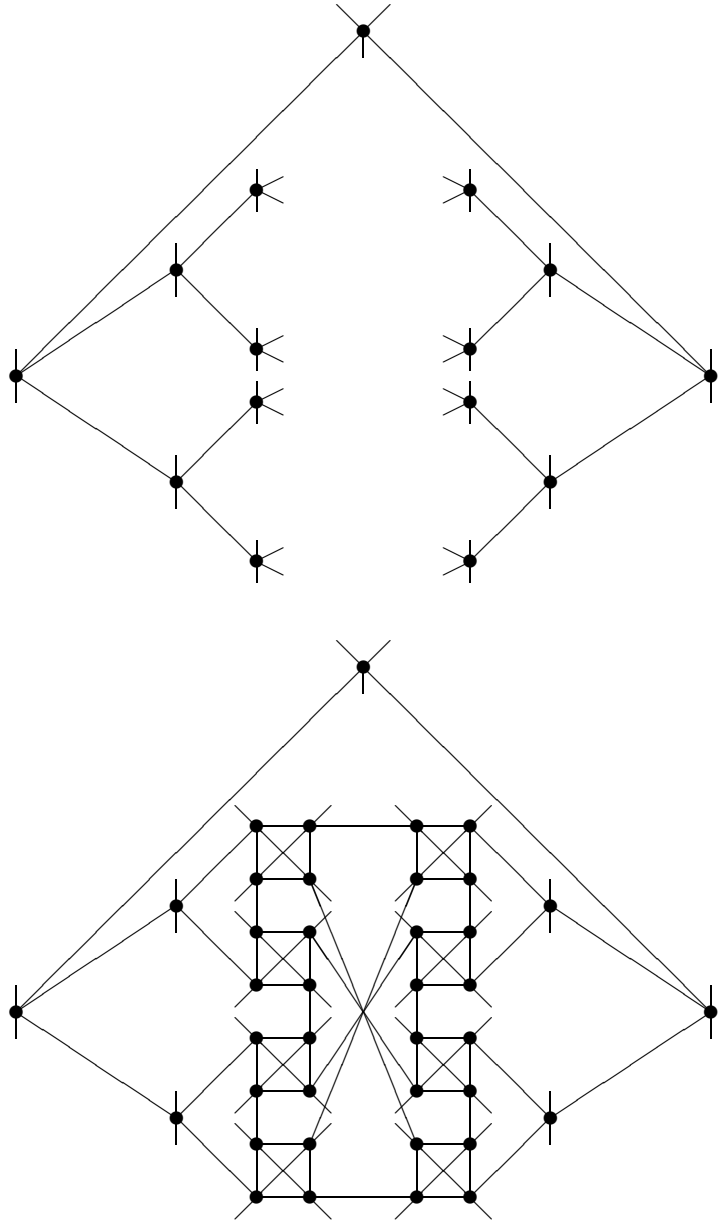

Figure 4: Graph $H_5(K_4)$ with 7860 vertices.

| $\Delta$ \ $D$ | 2 | 3 | 4 | 5 | 6 | 7 | 8 |
|---|---|---|---|---|---|---|---|
| 3 | $P$ 10 | $C_5 * F_4$ 20 | $vC$ 38 | $vC$ 70 | $GFS$ 130 | $CR^*$ 184 | $CR^*$ 320 |
| 4 | $K_3 * C_5$ 15 | $Allwr$ 41 | $C_5 * C_{19}$ 95 | $H'_3$ 364 | $H_3(K_3)$ **740** | $DH$ 1 155 | $DH^{**}$ 3 025 |
| 5 | $K_3 * X_8$ 24 | $Lente$ 70 | $Q_4(K_3)$ **186** | $H'_3 d$ 532 | $H_4(K_4)$ **2 754** | $DH$ 5 050 | $DH$ 13 780 |
| 6 | $K_4 * X_8$ 32 | $C_5 * C_{21}$ 105 | $DH^*$ 360 | $DH^*$ 1 200 | $H_5(K_4)$ **7 860** | $DH$ 18 205 | $DH$ 68 328 |
| 7 | $HS$ 50 | $DH^*$ 144 | $DH^*$ 600 | $DH$ 2 756 | $H_4(K_4) < H_5$ **10 566** | $DH$ 47 304 | $DH$ 165 228 |
| 8 | $P'_7$ 57 | $DH$ 220 | $DH$ 952 | $DH^*$ 4 704 | $H_7(K_6)$ **39 396** | $DH^{**}$ 110 889 | $DH$ 510 900 |
| 9 | $P'_8 d$ 74 | $Q'_8$ 585 | $Din$ 1 254 | $DH^*$ 7 260 | $H_8(K_6)$ **75 198** | $DH$ 218 130 | $DH^{**}$ 1 354 896 |
| 10 | $P'_9$ 91 | $Q'_8 d$ 650 | $DH$ 1 904 | $DH^*$ 12 288 | $H_9(K_6)$ **133 500** | $DH$ 504 710 | $2cy$ 3 000 000 |
| 11 | $P'_9 d$ 94 | $Q'_8 d$ 715 | $Q_7(T_4)$ 3 200 | $Din$ 16 578 | $H_7(T_4)$ 156 864 | $Din$ 914 414 | $Cam$ 4 773 696 |
| 12 | $P'_{11}$ 133 | $Q'_8 d$ 780 | $Q'_8 * X_8$ 4 680 | $Din$ 26 268 | $H_{11}(K_6)$ **355 812** | $Din$ 1 732 514 | $2cy$ 10 000 000 |
| 13 | $P'_{11} d$ 136 | $Q'_8 d$ 845 | $Q_9(T_4)$ 6 560 | $DH$ 36 290 | $H_9(T_4)$ 531 440 | $Cam$ 2 723 040 | $2cy$ 15 000 000 |
| 14 | $P'_{13}$ 183 | $Q'_8 d$ 910 | $Q_9(T_5)$ 8 200 | $DH$ 53 025 | $H_{13}(K_7)$ **806 636** | $K_1 \Sigma_8 H_{11}$ 6 200 460 | $Din$ 29 992 052 |

Table 1: Largest $(\Delta, D)$-graphs

**Graphs**

| | |
|---|---|
| $2cy$ | Connection of two cycles [?] |
| $Allwr$ | Special graphs found by Allwright [1] |
| $Cam$ | Cayley graphs found by Campbell [6] |
| $CR^*$ | Chordal rings found by Quisquater [16] |
| $vC$ | Compound graphs designed by von Conta [8] |
| $Din$ | Cayley graphs found by Dinneen [11] |
| $C_n$ | Cycle on $n$ vertices |
| $GFS$ | Special graph by Gómez, Fiol and Serra [14] |
| $H_q$ | Incidence graph of a regular generalized hexagon [2] |
| $HS$ | Hoffman-Singleton graph |
| $K_n$ | Complete graph |
| $Lente$ | Special graph designed by Lente, Univ. Paris Sud, France |
| $P$ | Petersen graph |
| $P_q$ | Incidence graph of projective plane [15] |
| $Q_q$ | Incidence graph of a regular generalized quadrangle [2] |
| $T$ | Tournament |

**Operations**

| | |
|---|---|
| $G * H$ | twisted product of graphs [3] |
| $Gd$ | duplication of some vertices of $G$ [10] |
| $B'$ | quotient of the bipartite graph $B$ by a polarity [9] |
| $B(K)$ | Substitution of vertices of a bipartite graph $B$ by complete graphs $K$ (this paper) |
| $B(T)$ | Compound using a bipartite graph $B$ and a tournament $T$ [13] |
| $B(K) < B$ | Compound of $B(K)$ and a bipartite graph $B$ [14] |
| $G\Sigma_i H$ | Various compounding operations [14] |
| $DH$ | Semidirect product of cyclic groups [12] |
| $DH*$ | Semidirect product of cyclic groups and $Z_n \times Z_n$ [12] |
| $DH**$ | Semidirect product $G \cdot G$ where G is a semidirect product of cyclic groups [12] |

Table 2: Graphs and operations in the table of large $(\Delta, D)$ graphs